\documentclass[12pt,letterpaper]{article}

\pdfoutput=1

\usepackage{graphics,epsfig,graphicx,color}
\graphicspath{{../Figures/}}
\usepackage[caption = false]{subfig}
\usepackage{float}
\usepackage{capt-of}

\usepackage{amssymb,latexsym}

\usepackage{setspace}

\usepackage{amsmath}

\usepackage{mathrsfs,amsfonts}

\usepackage{amsthm,amsxtra}

\usepackage[final]{showkeys}

\usepackage{stmaryrd}

\usepackage{algorithm}
\usepackage{algorithmicx}
\usepackage{algpseudocode}

\usepackage[openbib]{currvita}

\usepackage{enumitem}

\newif\ifPDF
\ifx\pdfoutput\undefined
	\PDFfalse
\else
	\ifnum\pdfoutput > 0
		\PDFtrue
	\else
		\PDFfalse
	\fi
\fi

\ifPDF
	\usepackage{pdftricks}
	\begin{psinputs}
		\usepackage{pstricks}
		\usepackage{pstcol}
		\usepackage{pst-plot}
		\usepackage{pst-tree}
		\usepackage{pst-eps}
		\usepackage{multido}
		\usepackage{pst-node}
		\usepackage{pst-eps}
	\end{psinputs}
\else
	\usepackage{pstricks}
\fi

\ifPDF
	\usepackage[debug,pdftex,colorlinks=true, 
	linkcolor=blue, bookmarksopen=false,
	plainpages=false,pdfpagelabels]{hyperref}
\else
	\usepackage[dvips]{hyperref}
\fi

\usepackage{fancyhdr}

\usepackage{comment}

\usepackage[toc,page]{appendix}

\usepackage{cancel}

\newcommand{\sgn}{\operatorname{sgn}}

\newcommand{\eps}{\varepsilon}

\newcommand{\bbR}{\mathbb R}

 \newcommand{\bn}{\mathbf n}

 \newcommand{\bx}{\mathbf x} 
\newcommand{\by}{\mathbf y} \newcommand{\bz}{\mathbf z}

\newcommand{\cG}{\mathcal G}

\newcommand{\angstrom}{\mbox{\normalfont\AA}}

\newenvironment{keywords}
{\noindent{\bf Key words.}\small}{\par\vspace{1ex}}
\newenvironment{AMS}
{\noindent{\bf AMS subject classifications 2010.}\small}{\par}

\title{An implicit boundary integral method for computing electric potential of macromolecules in solvent}

\author{
	Yimin Zhong\thanks{
		Department of Mathematics, 
		University of Texas, Austin, TX 78712; 
		yzhong@math.utexas.edu
	}
	\and 
	Kui Ren\thanks{
		Department of Mathematics and ICES,
		University of Texas, Austin, TX 78712; 
		ren@math.utexas.edu
	}
	\and 
	Richard Tsai\thanks{
		Department of Mathematics and ICES,
		University of Texas, Austin, TX 78712 and 
		KTH Royal Institute of Technology; 
		ytsai@math.utexas.edu
	}
}

\begin{document}

\maketitle



\begin{abstract}
A numerical method using implicit surface representations is proposed to solve the linearized Poisson-Boltzmann equation that {arises} in mathematical models for the electrostatics of molecules in solvent. The proposed method {uses} an implicit boundary integral formulation to {derive} a linear system defined on Cartesian nodes in a narrowband surrounding the closed surface that {separates} the molecule and the solvent. The needed implicit {surface} is constructed from the given atomic description of the molecules, by a sequence of standard level set algorithms. A fast multipole method is applied to accelerate the solution of the linear system. A few numerical studies involving some standard test cases are presented and compared to other existing results.
\end{abstract}






\section{Introduction}
\label{SEC:Intro}

The mathematical modeling and numerical simulation of electrostatics of charged macromolecule-solvent {systems} have been extensively studied in recent years, due to their importance in many branches of electrochemistry; see, for instance, ~\cite{BaFa-Book00,DaMc-CR90,FeBr-COSB04,FeOnLeImCaBr-JCC04,HoNi-Science95,LuZhHoMc-CiCP08,Mackerell-JCC04,MaCrTr-JPCB09,OrLu-CR00,RoSi-BC99,ToMeCa-CR05,ZhGaFrLe-JCC01} and references therein for recent overviews of the developments in the subject.

There are roughly two classes of mathematical models for such macromolecule-solvent {systems}, depending on how the effect of the solvent is modeled: explicit solvent models in which solvent molecules are treated explicitly, and implicit solvent models in which the solvent is represented as a continuous medium. While explicit solvent models are believed to be more accurate, they are computationally intractable when modeling large systems. Implicit models are therefore often an alternative for large simulations, see~\cite{Baker-COSB05,ChBrKh-COSB08,CrTr-CR99,LaHe-JCP11,ZhWiAl-PB11} and references therein for recent advances. The Poisson-Boltzmann model is one of the popular implicit solvent models in which the solvent is treated as a continuous high-dielectric medium~\cite{CaWaZhLu-JCP09,Chipman-JCP04,FoBrMo-JMR02,ImBeRo-CPC98,LiXi-CMS15,LuDaGi-JCC02,ReChThScMaZhBa-QRB10,RoAlHo-JPCB01,WaTaTaLuLu-CiCP08}. This model, and many variants of it, has important applications, for instance in studying biomolecule dynamics of large proteins~\cite{BaSeJoHoMc-PNAS98,Baker-ME04,Bardhan-JCP11,BoAnOr-PRL97,GrPeVa-JCP07,ZhWiAl-PB11}. Many efficient and accurate computational schemes for the numerical solution of the model have been developed~\cite{BaHoWa-JCC00,BaScDu-PRE09,Bardhan-JCP09,ChChChGeWe-JCC11,DiSp-JCP16,Geng-JCP13,GeKr-JCP13,HeGi-JCP11,LiXi-CMS15,LuChHuMc-PNAS06,WeRuHi-JCP10,ZhLuChHuPiSuMc-CiCP13}.

\begin{figure}
\begin{centering}
\includegraphics[scale = 0.35]{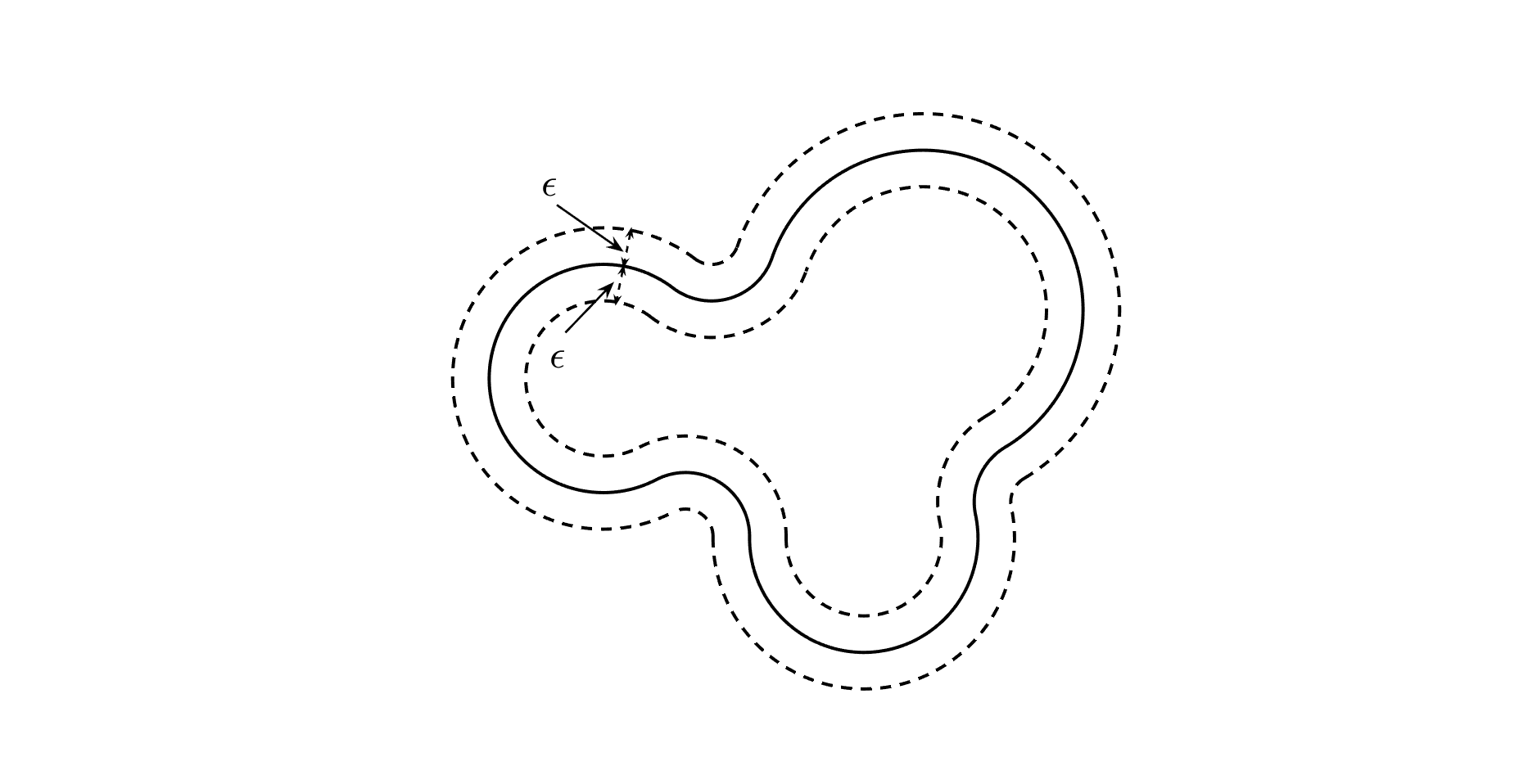}
\par\end{centering}
\caption{A view of the ``\emph{solvent excluded surface}'' in 2D is shown by the middle solid curve, and the narrowband  $\Gamma_\epsilon$ is shown here by the space bounded between the dashed curves. }
\label{fig:2aid_illu}
\end{figure}

To introduce the Poisson-Boltzmann model, let us assume that the macromolecule has $N_c$ atoms centered at $\{\bz_j\}_{j=1}^{N_c}$, with radii $\{r_j\}_{j=1}^{N_c}$ and charge number $\{q_j\}_{j=1}^{N_c}$ respectively. Let $\Gamma$ be the closed surface that separates the region occupied by the macromolecule and the rest of the space. The typical choice of $\Gamma$ is the so-called \emph{solvent excluded surface}, which is defined as the boundary of the region outside the macromolecule which is accessible by a probe sphere with some small radius, say $\rho_0$; see Figure~\ref{fig:2aid_illu} for an illustration. We use $\Omega$ to denote the region surrounded by $\Gamma$ that includes the macromolecule. 

We use a single function $\psi$ to denote the electric potential inside and outside of $\Omega$. In the Poisson-Boltzmann model, $\psi$ solves the Poisson's equation for point {charges} inside $\Omega$, that is,
\begin{equation*}
	-\nabla\cdot (\eps_I\nabla \psi(\bx)) = \sum_{k=1}^{N_c} q_k \delta(\bx - \bz_k),\ \ \mbox{in}\ \ \Omega
\end{equation*}
where $\eps_I$ denotes the dielectric constant in $\Omega$. Outside $\overline\Omega$, that is in the solvent that {excludes} the interface $\Gamma$, $\psi$ solves the Poisson's equation for a continuous distribution of charges that models the effect of the solvent, that is,
\begin{equation*}
-\nabla\cdot\left(\epsilon_E\nabla\psi\right)=  \rho_B(T, \bx,\psi(\bx)),\ \ \mbox{in}\ \ \bbR^3\setminus\overline{\Omega}
\end{equation*}
where $\eps_E$ denotes the dielectric constant of the solvent, which often has much higher value than that of the macromolecule, $\eps_E\gg \eps_I$. The source term $\rho_B$ is a nonlinear function coming from the Boltzmann distribution with $T$ denoting the temperature of the system. More precisely, for solvent containing $m$ ionic species,
\[
\rho_B(T, \bx,\psi(\bx)):= e_c \sum_{i=1}^m   c_i \bar q_i e^{-e_c \bar q_i\psi(\bx)/k_BT},~~~\bx\in\bbR^3\setminus\overline{\Omega}
\]
where $c_i,\bar q_i$ are the concentration and charge of the $i$th ionic species, $e_c$ is the electron charge, $k_B$ is the Boltzmann constant, and $T$ is the absolute temperature. 

The nonlinear term $\rho_B(T,\bx,\psi)$ in the Poisson-Boltzmann system poses significant challenges in the computational solution of the system. In many practical applications, it is replaced by the linear function $-{\bar\kappa_T}^2\psi(\bx)$ where the parameter $\bar{\kappa}_T =\sqrt{ \frac{2 e^2 \mathbb{I}}{k_B T}}$ is called the Debye-Huckel screening parameter with $k_B$, $e$, and $\mathbb{I}$ being the Boltzmann constant, the unit charge, and the ionic strength respectively. This leads to the linearized Poisson-Boltzmann equation (PBE) for the electrostatic potential $\psi$. It takes the following form
\begin{equation}\label{EQ:PB}
\begin{aligned}
	-\nabla\cdot (\eps_I\nabla \psi(\bx)) &= \sum_{k=1}^{N_c} q_k \delta(\bx - \bz_k),&&\mbox{in}\ \ \Omega, \\
	-\nabla\cdot (\eps_E \nabla\psi(\bx)) &= -\bar{\kappa}_T^2 \psi(\bx), &&\mbox{in}\ \ {\overline{\Omega}}^c,\\
	\psi(\bx)_{|\Gamma_+} &= \psi(\bx)_{|\Gamma_-} , &&\mbox{on}\ \ \Gamma,\\
	\eps_E{\frac{\partial \psi}{\partial n}}_{|\Gamma_+} &= \eps_I{\frac{\partial \psi}{\partial n}}_{|\Gamma_-},&& \mbox{on}\ \ \Gamma,\\
	|\bx|\psi(\bx) \to 0, & \ \ \ \ |\bx|^2|\nabla\psi(\bx)|\to 0, &  & \mbox{as}\ |\bx|\to \infty.
	\end{aligned}
\end{equation}
Here the operator $\partial/\partial n\equiv \bn(\bx)\cdot\nabla$ denotes the usual partial derivative at $\bx\in\Gamma$ in the outward normal direction $\bn(\bx)$ (pointing from {$\Gamma$} outward). The usual continuity conditions, continuity of the potential and the flux across $\Gamma$, are assumed, and the radiation condition, which requires $\psi$ decay to zero far away from the macromolecule, is needed to ensure the uniqueness of solutions to the {linearized} Poisson-Boltzmann equation. See e.g. ~\cite{BaHoWa-JCC00,BaScDu-PRE09,Bardhan-JCP11,CaWaZhLu-JCP09,Geng-JCP13,GeKr-JCP13,LaHe-JCP11,LuZhHoMc-CiCP08,RoAlHo-JPCB01,WaTaTaLuLu-CiCP08,WeRuHi-JCP10}.

Computational solution of the {linearized} Poisson-Boltzmann equation~\eqref{EQ:PB} in practically relevant configurations turns out to be quite challenging. Different types of numerical methods, including for instance finite difference methods~\cite{Baker-COSB05,BrBrOlStSwKa-JCC83,ImBeRo-CPC98,GeYuWe-JCP07,MaBrWaDaLuIlAnGiBaScMc-CPC95,NiBhHo-BJ93}, finite element methods~\cite{BaHoWa-JCC00,HoSa-JCC95,XiJi-JCP16,XiJiBrSc-SIAM12,XiYiXi-JCC17,YiXi-JCP15}, boundary element methods~\cite{AlBaWhTi-JCC09,BaChRa-SIAM11,BoFeZh-JPCB02,JuBovavaBe-JCP91,LiSu-BJ97,LuChHuMc-PNAS06,LuChMc-JCP07,ZaMo-JCC90}, and many more hybrid or specialized methods~\cite{BoFe-JCC04,VoGrSc-ACS92} have been developed; see~\cite{LuZhHoMc-CiCP08} for the recent survey on the subject. Each method has its own advantages and disadvantages. Finite difference methods are easy to implement. They are the methods used in many existing software packages~\cite{BrBrOlStSwKa-JCC83,ImBeRo-CPC98,GeYuWe-JCP07,MaBrWaDaLuIlAnGiBaScMc-CPC95,NiBhHo-BJ93}. 
Finite difference methods, except the ones that are based on adaptive oct-tree structures~\cite{HeGi-JCP11, MIRZADEH20112125, MiThHeBoGi-CiCP13}, use uniform Cartesian grids and require special care for implementing the interface conditions to high order while maintaining stability.
Finite element methods provide more flexibility with the geometry. However, like the finite difference methods, they often suffer from issues such as large memory storage requirement and low solution speed when dealing with large problems. Moreover, both finite difference and finite element methods need to truncate the domain in some way, therefore the radiation condition is not satisfied exactly. Boundary element methods are based on integral formulations of the Poisson-Boltzmann equation. They require only the discretization of the solvent excluded surface, i.e. $\Gamma$, not the macromolecule and solvent domains. The radiation condition is usually exactly, although implicitly, integrated into the integral form to be solved. However, the matrix systems {resulting} from boundary element formulations are often dense. Efficient acceleration, for instance preconditioning, techniques are needed to accelerate the solution of such dense systems.


In this work, we propose a fast numerical {method} for solving the interface/boundary value problem of the linearized Poisson-Boltzmann equation~\eqref{EQ:PB}. The method is derived from the implicit boundary integral formulation~\cite{KTT} of~\eqref{EQ:PB} and relies on some of the classical level set algorithms~\cite{LevelSet_OsherFedkiw,osher_sethian88} for computing the implicit interfaces and the needed geometrical information. All the involved computational procedures are defined on an underlying uniform Cartesian grid. Thus the proposed method {inherits most} of the flexibilities of a level set algorithm.
On the other hand, since the method is derived from a boundary integral formulation of~\eqref{EQ:PB}, it treats the interface conditions and far field conditions in a less involved fashion compared to the standard level set algorithms for similar problems. As such {types} of implicit boundary integral approaches are relatively new, we describe in detail how to set up a linear system and where a fast multipole method can be used for acceleration of the common matrix-vector multiplications in the resulting linear system. We demonstrate in our simulations involving non-trivial molecules defined by tens of thousands atoms that standard ``kernel-independent'' fast multipole methods~\cite{fong2009black} can be used easily and effectively as in a standard boundary integral method.




{We conclude the introduction with the following remarks. The linearized Poisson-Boltzmann model provides a sufficiently accurate approximation in many cases, in particular when the solution's ionic strength is relatively low; see \cite{FoBrMo-JMR02} and references therein. In cases where the linearized model is not accurate enough, solving the nonlinear Poisson-Boltzmann is necessary. With appropriate far field conditions, finite difference or finite element methods provide a way to compute solutions in such case; see~\cite{MiThHeBoGi-CiCP13} and references therein. The computational method we develop in this work can potentially be combined with an iterative scheme for nonlinear equations, such as methods of Newton's type~\cite{Kelley-Book95}, to solve the nonlinear Poisson-Boltzmann equation. At each iteration of the nonlinear solver, the proposed method can be adapted to solve the linearized problem as long as the coefficients involved, mainly the dielectric coefficients, are constants as currently assumed in our algorithm. }

The rest of this paper is organized as follows. We first introduce in Section~\ref{SEC:IBIM} the implicit boundary integral formulation of the linearized Poisson-Boltzmann system~\eqref{EQ:PB}. We then present the details of the implementation of the method in Section~\ref{SEC:Num Impl}. In Section~\ref{SEC:Num Exp}, we present some numerical simulation results to demonstrate the performance of the algorithm. Concluding remarks are then offered in Section~\ref{SEC:Concl}.

\section{The implicit boundary integral formulation}
\label{SEC:IBIM}

The numerical method we develop in this work is based on a boundary integral formulation of the linearized Poisson-Boltzmann equation that is developed in~\cite{JuBovavaBe-JCP91}.

\subsection{Boundary integral formulation} 

Throughout the rest of the paper, all the coefficients involved in the equations are assumed to be constant, i.e. independent of the spatial variable. We define $\kappa = \bar\kappa_T/\sqrt{\epsilon_E}$, and introduce the fundamental solutions
\begin{equation*}
	G_0(\bx, \by) = \frac{1}{4\pi |\bx - \by|} \ \ \mbox{and}\ \ G_{ \kappa}(\bx, \by) = \frac{e^{-\kappa|\bx-\by|}}{4\pi|\bx - \by|}
\end{equation*}
to the Laplace equation and the one with the linear lower order term $-\bar\kappa^2_T\psi$  in \eqref{EQ:PB}.

Following the standard way of deriving boundary integral equations, we apply Green's theorem to the system formed by (i) the first equation in ~\eqref{EQ:PB} and the equation for $G_0$, and (ii) the second equation in~\eqref{EQ:PB} and the equation for $G_\kappa$, taking into account the interface and the radiation conditions. A careful routine calculation leads to the following boundary integral equations for the potential $\psi$ and its normal derivative $\psi_n\equiv \partial\psi/\partial n$ on $\Gamma$:
\begin{equation}\label{EQ:PB Intg0}
\begin{aligned}
	\frac{1}{2}\psi(\bx) + \int_{\Gamma} \left(\frac{\partial G_0(\bx, \by)}{\partial n(\by)}\psi(\by)  -  G_0(\bx, \by)\psi_n(\by)\right) d\by &= \sum_{k=1}^{N_c} \frac{q_k}{\epsilon_I} G_0(\bx, \bz_k),\\
	\frac{1}{2}\psi(\bx) - \int_{\Gamma} \left(\frac{\partial G_{ \kappa}(\bx, \by)}{\partial n(\by)}\psi(\by) -\frac{\epsilon_I}{\epsilon_E} G_{ \kappa}(\bx, \by)\psi_n(\by)\right) d\by &= 0.
\end{aligned}
\end{equation}
This system of boundary integral equations is the starting point of many existing numerical algorithms for the linearized Poisson-Boltzmann equation. 

In our algorithm, we adopt the integral formulation proposed in~\cite{JuBovavaBe-JCP91}. This formulation reads:
\begin{equation}\label{EQ:PB Intg1}
\begin{aligned}
	\frac{1}{2}\left(1 + \frac{\epsilon_E}{\epsilon_I}\right)\psi(\bx) &+ \int_{\Gamma} \left(\frac{\partial G_0(\bx, \by)}{\partial n(\by)} - \frac{\epsilon_E}{\epsilon_I}\frac{\partial G_{\kappa}(\bx, \by)}{\partial n(\by)}\right)\psi(\by) d\by \\&- \int_{\Gamma} \left(G_0(\bx, \by) - G_{\kappa}(\bx, \by)\right)\psi_n(\by) d\by = \sum_{k=1}^{N_c} \frac{q_k}{\epsilon_I}G_0(\bx, \bz_k) ,
\end{aligned}
\end{equation}
\begin{equation}\label{EQ:PB Intg2}
\begin{aligned}
	\frac{1}{2}\left(1+\frac{\epsilon_I}{\epsilon_E}\right)\psi_n(\bx) &+ \int_{\Gamma} \left(\frac{\partial^2 G_0(\bx,\by)}{\partial n(\bx)\partial n(\by)} -\frac{\partial^2 G_{\kappa}(\bx,\by)}{\partial n(\bx)\partial n(\by)} \right)\psi(\by) d\by \\ &-\int_{\Gamma} \left(\frac{\partial G_0(\bx, \by)}{\partial n(\bx)} - \frac{\epsilon_I}{\epsilon_E}\frac{\partial G_{\kappa}(\bx,\by)}{\partial n(\bx)}\right)\psi_n(\by) d\by = \sum_{k=1}^{N_c}\frac{q_k}{\epsilon_I}\frac{\partial G_0(\bx, \bz_k)}{\partial n(\bx)}.
\end{aligned}
\end{equation}
The first equation in this formulation,~\eqref{EQ:PB Intg1}, is simply the linear combination of the two equations in~\eqref{EQ:PB Intg0}, while the second equation in this formulation,~\eqref{EQ:PB Intg2}, is nothing but the linear combination of the derivatives of the two equations in~\eqref{EQ:PB Intg0}. It is shown in~\cite{JuBovavaBe-JCP91} that the potentially hypersingular integral in~\eqref{EQ:PB Intg2}, involving the second derivatives of $G_0$ and $G_\kappa$ is actually integrable on $\Gamma$, thanks to the fact that
\[
	\frac{\partial^2 G_0(\bx,\by)}{\partial n(\bx)\partial n(\by)} -\frac{\partial^2 G_{\kappa}(\bx,\by)}{\partial n(\bx)\partial n(\by)}\sim O(|\bx-\by|^{-1}),\ \ |\bx-\by|\rightarrow 0.
\]
Moreover, when $\kappa=0$,~\eqref{EQ:PB Intg1} is decoupled from~\eqref{EQ:PB Intg2}, and the latter provides an explicit formula for evaluating $\partial\psi/\partial n$ using $\psi$.

The main benefit of the formulation~\eqref{EQ:PB Intg1}-\eqref{EQ:PB Intg2} is that it typically leads to, after discretization, linear systems with smaller condition numbers than the formulation in~\eqref{EQ:PB Intg0}. The typical boundary element methods for this system (and others) require careful triangulation of the interface $\Gamma$; see e.g. ~\cite{AlBaWhTi-JCC09,BaChRa-SIAM11,BoFeZh-JPCB02,JuBovavaBe-JCP91,LiSu-BJ97,LuChMc-JCP07,ZaMo-JCC90}. In the next subsection, we describe our method to discretize the boundary integral system~\eqref{EQ:PB Intg1} and \eqref{EQ:PB Intg2} on a subset of a uniform Cartesian grid nodes in a narrowband surrounding $\Gamma$, without the need to parameterize $\Gamma$.

	
\subsection{Implicit boundary integral method}
\label{sec:IBIM}

Let the interface $\Gamma$ be a closed, $C^{1,\alpha}$ surface (in two or three {dimensions})  with $\alpha>0$  so that the distance function to $\Gamma$ is differentiable in a neighborhood around it. Let $d_\Gamma$ denote the signed distance function to $\Gamma$ that takes the negative sign for points inside the region enclosed by $\Gamma$, and  $\Gamma_\epsilon$ denote the set of points whose distance to $\Gamma$ is smaller than $\epsilon$.
An implicit boundary integral formulation of a surface integral defined on $\Gamma$ is derived by
projecting points in $\Gamma_\epsilon$ onto their closest points on $\Gamma$. With the distance function to $\Gamma$, the projection operator can be evaluated by
\begin{equation}
P_\Gamma \bx := \bx - d_\Gamma(\bx)\nabla d_\Gamma(\bx).
\end{equation}
When $\epsilon$ is smaller than the maximum principal curvatures of $\Gamma$, the closest point projection is well-defined in $\Gamma_\epsilon$. 
		
An implicit boundary integral method (IBIM) \cite{KTT} is built upon the following identity:
\begin{equation}\label{eq:IBIM-equivalence}
	I_\Gamma[f]:=\int_{\Gamma} f(\bx) ds(\bx) = \int_{\Gamma_{\eps}} f(P_\Gamma \bx) \delta_{\eps}(d_\Gamma(\bx))J(\bx) d\bx,
\end{equation}
which reveals the equivalence between the surface integral and its extension into a volume integral. 
We shall call the integral over $\Gamma_{\eps}$ an implicit boundary integral.
In this implicit boundary integral, one has
\begin{enumerate}
\item The extension of $f(\bx)$ as a constant along the normal of $\Gamma$ at $\bx$. 
\item The Jacobian $J(\bx)$ which accounts for the change of variables between $\Gamma$ and the {level set surface of $d_{\Gamma}$} that passes through $\bx$.
\item A weight function, $\delta_\epsilon$ compactly supported on $[-\eps, \eps]$ satisfying
\begin{equation}
	\int_{-\eps}^{\eps} \delta_{\eps}(\eta) d\eta = 1.
\end{equation} 
\end{enumerate}

In $\mathbb{R}^3$, The Jacobian $J$ takes the explicit {form}
\begin{equation}
	J(\bx) = 1 - d_\Gamma(\bx)\Delta d_\Gamma(\bx) + d_\Gamma(\bx)^2 \langle \nabla d, \nabla^2 d_\Gamma \nabla d_\Gamma\rangle.
\end{equation}
It can be further related to the products of the singular values of the Jacobian matrix of $P_\Gamma$, which provides an alternative, and in some cases easier way, for the computation of $J$. See~\cite{kublik2016integration}.

In the application of  interest, the distance function to a solvent excluded surface will be twice differentiable almost everywhere in some narrowband around the surface, together with setting  $J\equiv 1$, the proposed method is well-defined in there. 
In fact, the requirement on the regularity of the surface (and its signed distance function) 
can be further relaxed if the weight function is an even function  and possesses enough vanishing moments.  
It is shown in~\cite{kublik2016extrapolative} that with $J\equiv 1$, the cosine weight function \eqref{eq:cosine-kernel} and
$C^1$ integrands which are not necessarily constant along surface normals, $I_\Gamma[f]$ is approximated to second order 
if $\epsilon$.
It is further shown in~\cite{kublik2016extrapolative} that if the weight function has more than two vanishing moments, one may replace the Jacobian by $J\equiv 1$ while keeping the equality in~\eqref{eq:IBIM-equivalence} valid, even for piecewise smooth surfaces containing corners and creases.
	
Numerically we approximate the implicit boundary integral by embedding the computational domain $\Omega$ into the rectangle $U = [a,b]^n$, and subdivide $U$ into the uniform grid $U_h=h\mathbb{Z}^n\cap [a,b]^n$ with grid size $h = (b-a)/N$ along each coordinate direction and $\bx_i$ at each grid point. We approximate the integral by  
\begin{equation}\label{EQ:IBIM sum}
	I_\Gamma[f] \approx S^h_{\Gamma_\epsilon} [f] :=\sum_{\bx_i\in \Gamma_{\eps}} f(\bx_i^{\ast}) \delta_{\eps}(d_\Gamma(\bx_i)) J(\bx_i) h^n
\end{equation}
where $\bx_i^{\ast} = \bx_i - d_\Gamma(\bx_i) \nabla d_\Gamma(\bx_i)$ is the projection of $\bx_i$ onto $\Gamma$.

Thus, a typical second kind integral equation of the form
\begin{equation}\label{eq:generic-second-kind-BIE}
	g(\bx) = \lambda \beta(\bx) + \int_\Gamma K(\bx,\by)\beta(\by) ds(\by),~~~\bx\in\Gamma,
\end{equation}
can be approximated on $U_h$ using the IBIM formulation. One would derive a linear system for the unknown function $\bar\beta$ defined on the grid nodes in $\Gamma_\epsilon$:
\begin{equation}\label{eq:IBIM-linear-system}
	g(P_\Gamma \bx_i) = \lambda \bar\beta(\bx_i) + h^n \sum_{\by_j\in\Gamma_\epsilon\cap U_h} K(P_\Gamma\bx_i,P_\Gamma \by_j)\bar\beta(\by) \delta_\epsilon( d_\Gamma(\by_j) ) J(\by_j),~~~\bx_i\in\Gamma_\epsilon\cap U_h,
\end{equation}
with the property that as $h\rightarrow 0$
\[
	\bar\beta( \bx_i ) \longrightarrow \beta( P_\Gamma \bx_i ),~~~\forall \bx_i\in\Gamma_\epsilon\cap U_h;
\]
i.e. the solution to the linear system \eqref{eq:IBIM-linear-system} converges to ``the function which is the constant extension along the surface normal'' of the solution of~\eqref{eq:generic-second-kind-BIE}; see more discussions in~\cite{chen2016-Helmholtz,KTT}.

In the context of this paper, equations~\eqref{EQ:PB Intg1} and~\eqref{EQ:PB Intg2} will be discretized into
\begin{equation}\label{eq:IBIM-PBE}
\begin{aligned}
\frac{1}{2}\lambda_{1}\bar\psi_{~}(\bx_{i})+h^{3}\sum_{j}K_{11}(\bx_{i},\by_{j})\omega_j\bar\psi(y_{j})-h^{3}\sum_{j}K_{12}(\bx_{i},\by_{j})\omega_j{\bar\psi}_{n}(\by_{j})	&=g_{1}(\bx_{i}),\\
\frac{1}{2}\lambda_{2}\bar\psi_{n}(\bx_{i})+h^{3}\sum_{j}K_{21}(\bx_{i},\by_{j})\omega_j\bar\psi(y_{j})-h^{3}\sum_{j}K_{22}(\bx_{i},\by_{j})\omega_j{\bar\psi}_{n}(\by_{j})	&=g_{2}(\bx_{i}),
\end{aligned}
\end{equation}
where 
$$\lambda_1=\frac{1}{2}\left(1 + \frac{\epsilon_E}{\epsilon_I}\right), \lambda_2=\frac{1}{2}\left(1 + \frac{\epsilon_I}{\epsilon_E}\right),$$ 
\begin{equation}\label{eq:weights}
\omega_j:=J(\by_j)\delta_\epsilon(d_\Gamma(\by_j)),
\end{equation}  
\[
g_1(\bx_i) :=\sum_{k=1}^{N_c} \frac{q_k}{\epsilon_I}G_0(\bx_i, \bz_k),
g_2(\bx_i) := \sum_{k=1}^{N_c}\frac{q_k}{\epsilon_I}\frac{\partial G_0(\bx_i, \bz_k)}{\partial n(\bx_i)},
\]
 and $K_{11}, K_{12}, K_{21}, K_{22}$ are respectively the regularized versions of the following weakly singular kernels:
$$\frac{\partial G_0(P_\Gamma\bx, P_\Gamma\by)}{\partial n(P_\Gamma\by)} - \frac{\epsilon_E}{\epsilon_I}\frac{\partial G_{\kappa}(P_\Gamma\bx, P_\Gamma\by)}{\partial n(P_\Gamma\by)}, 
G_0(P_\Gamma\bx, P_\Gamma\by) - G_{\kappa}(P_\Gamma\bx, P_\Gamma\by),$$ 
$$\frac{\partial^2 G_0(P_\Gamma\bx,P_\Gamma\by)}{\partial n(P_\Gamma\bx)\partial n(P_\Gamma\by)} -\frac{\partial^2 G_{\kappa}(P_\Gamma\bx,P_\Gamma\by)}{\partial n(P_\Gamma\bx)\partial n(P_\Gamma\by)},~ 
\textrm{and}~ 
\frac{\partial G_0(P_\Gamma\bx, P_\Gamma\by)}{\partial n(P_\Gamma\bx)} - \frac{\epsilon_I}{\epsilon_E}\frac{\partial G_{\kappa}(P_\Gamma\bx,P_\Gamma\by)}{\partial n(P_\Gamma\bx)},$$
A simple regularization that we used in our numerical implementation is described below in the next subsection.   


This formulation provides a convenient computational approach for computing boundary integrals, where the boundary is naturally defined implicitly, as a level set of a continuous function, 
and is difficult to {parameterize}. The geometrical information about the boundary is restricted to the computation of the Jacobian $J$ and the closest point extension of the integrand $f$ --- both of which can be approximated easily by simple finite differencing applied to the distance function $d_\Gamma(\bx)$ at grid point $\bx_i$ within $\Gamma_{\eps}$.
Furthermore, the smoothness of the weight function $\delta_\epsilon$, 
along with the smoothness of the {integrand} will allow for higher order in $h$ approximation of $I[f]$ by simple Riemann sum $S^h_{\Gamma_\epsilon}[f]$, see for example the discussion in~\cite{kublik2016integration}.

\subsubsection{Regularization of the kernels}\label{sec:regularization}

While all the kernels (the Green's functions and the particular linear combinations of them) that {appear} in \eqref{EQ:PB Intg1}-\eqref{EQ:PB Intg2} are formally integrable, an additional treatment for the singularities is needed in the numerical computation when $\bx$ and $\by$ are close.
Typically, the additional treatment corresponds to either local change of variables so that in the new variables the singularities do not exist or mesh refinement for control of numerical error amplification (particularly for Nystr\"om methods). 
The proposed simple discretization of the Implicit Boundary Integral formulation on uniform Cartesian grid can be viewed as an extreme case of Nystr\"om method, in which no mesh refinement is involved (and thus no control of numerical errors if the singularities of the kernels are left untreated). Therefore we need to {regularize} the kernels analytically and locally only when $P_\Gamma \bx$ and $P_\Gamma \by$ are sufficiently close with respect to the grid spacing. 

In the following, for brevity of the displayed formulas, let $\bx^*:=P_\Gamma \bx$, $\by^*:=P_\Gamma \by$ and 
\[
K_\theta(\bx,\by):=\frac{\partial G_0(\bx^*, \by^*)}{\partial n(\by^*)} - \theta\frac{\partial G_{\kappa}(\bx^*, \by^*)}{\partial n(\by^*)},~~~\theta\in\mathbb{R}.
\]
The regularization that we will use involves a small parameter $\tau>0$ and is defined by
\begin{equation}\label{eq:reg-kernel}
	K_{\theta}^{reg}(\bx, \by)= 
	\begin{cases}
	\overline{K_\theta}(\bx),  & \text{if } \|\bx^* - \by^*\|_{P} < \tau,\\
	K_\theta(\bx, \by),              & \text{otherwise},
	\end{cases}
\end{equation}
where $\|\bx^* - \by^*\|_P$ is the distance between projections of $\bx^*$ and $\by^*$ onto the tangent plane at $\bx^*$. 
$\overline{K_\theta}(\bx)$ is the average of $K_\theta(\bx,\cdot)$ defined as
\begin{equation}
	\overline{K_\theta (\bx)}= \frac{1}{V(\bx^*; \tau)}\int_{V(\bx;\tau)} K_\theta(\bx, \bz) ds(\bz), 
\end{equation}
where $V(\bx^*;\tau)$ is the disc of radius $\tau$ in the tangent plane of $\Gamma$ at $\bx^*$.

Thus, 
\begin{equation}
	K_{11}(\bx,\by):= K_\theta^{reg}(\bx,\by),~~~\textrm{with}~~~\theta=\frac{\epsilon_E}{\epsilon_I},~~\overline{K_\theta}(\bx)=0,
\end{equation}
\begin{equation}
	K_{22}(\bx,\by):= K_\theta^{reg}(\bx,\by),~~~\textrm{with}~~~\theta=\frac{\epsilon_I}{\epsilon_E},~~\overline{K_\theta}(\bx)=0, 
\end{equation}
Similarly, the averages of $G_0-G_\kappa$ and $\frac{\partial^2 G_0}{\partial n(\bx^*)\partial n(\by^*)} -\frac{\partial^2 G_{\kappa}}{\partial n(\bx^*)\partial n(\by^*)}$ are computed and we define: 
\begin{equation}\label{eq:reg-kernel-a}
	K_{12}(\bx, \by)= 
	\begin{cases}
	\frac{e^{-\kappa \tau} - 1 +\kappa\tau}{2\pi \kappa \tau^2},  & \text{if } \|\bx^* - \by^*\|_{P} < \tau,\\
	G_0(\bx^*, \by^*) - G_{\kappa}(\bx^*, \by^*),              & \text{otherwise},
	\end{cases}
\end{equation}

\begin{equation}\label{eq:reg-kernel-b}
	K_{21}(\bx, \by)= 
	\begin{cases}
	0,  & \text{if } \|\bx^* - \by^*\|_{P} < \tau,\\
	\frac{\partial^2 G_0(\bx^*,\by^*)}{\partial n(\bx^*)\partial n(\by^*)} -\frac{\partial^2 G_{\kappa}(\bx^*,\by^*)}{\partial n(\bx^*)\partial n(\by^*)},              & \text{otherwise}.
	\end{cases}
\end{equation}
Finally, we refer the readers to \cite{chen2016-Helmholtz} for a recent approach for dealing with hypersingular integrals via extrapolation.


\section{The proposed algorithm}
\label{SEC:Num Impl}

\begin{figure}[!htb]
\centering
\includegraphics[clip, trim={0.cm 5.5cm 0cm 13cm}, width=0.70\textwidth]{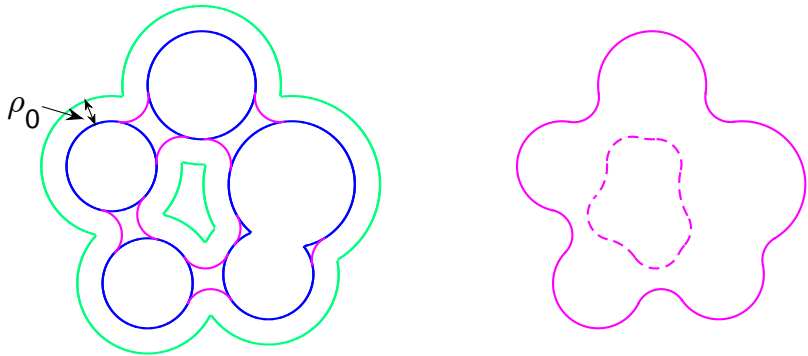}
\caption{The construction of the ``solvent excluded surface'' (SES) of a fictitious molecule defined by five atoms. The final SES is shown as the solid curve on the right plot. The  van der Waals surface corresponding to the molecule is shown by the blue curve. The green curve is the ``solvent accessible surface'', from which an inward eikonal flow will shrink it by a distance of $\rho_0$ to {arrive} at the pink curves (solid and dashed). 
The dashed pink curve on the right plot shows that boundary of the cavity enclosed by the molecules. It is removed from our computation.}
\label{fig:surfs}
\end{figure}

The proposed algorithm consists of a few stages which are outlined below:
\begin{enumerate}[label= Stage (\arabic*),leftmargin=3\parindent]\label{enum:stages}
\item {\bf Preparation of the signed distance function} to the ``solvent excluded surface'' on a uniform Cartesian grid.

This {includes definition of an initial level} set function (Section~\ref{sec:vdW-F}), followed by an ``inward''  eikonal flow of the level set function (Section~\ref{sec:inward-eikonal-flows}). After the eikonal flow, we apply a step that removes from the implicit surface the interior cavities which are not accessible to solvent (Section~\ref{sec:cavity-removal}).
See Figure~\ref{fig:surfs} for an illustration of this process and the various surfaces involved in it.

Finally we apply the reinitialization procedure (Section \ref{sec:reinit}) to the level set function obtained from cavity removal.
At the end of this stage, one shall obtain the signed distance function to the ``solvent excluded surface" on which the linearized Poisson-Boltzmann boundary integral equation (BIE) is solved. The constructed signed distance function has the same sign as the function $F$, defined in \eqref{vdW-levelset-fn}   that is used to defined the van der Waals surface, at the prescribed molecule centers.

\item {\bf Preparation of the linear system.} 

This involves computation of geometrical information, including the closest point mapping and  the Jacobian (Section~\ref{sec:weight-computations}). 

\item {\bf Solution of the linear system} via \texttt{GMRES} with a fast multipole acceleration for the matrix-vector multiplication (Section~\ref{sec:weight-computations}).

At the end of this stage, one {obtains} the density $\bar\psi$ defined on  the grid nodes lying in $\Gamma_\epsilon$. This density function will be used in the evaluation of the polarization energy.

\item{\bf Evaluation of surface area and polarization energy.} (Sections~\ref{sec:surface-area}~and~\ref {sec:polarization-energy-computation})

The surface area is computed by using $f\equiv1$ in~\eqref{eq:IBIM-equivalence} and the polarization energy is computed through the density $\bar{\psi}$ by the implicit boundary integral method.
\end{enumerate}

All computations will be performed on functions defined on $U_h$. The inward eikonal flow and the reinitialization {in Stage (1)} 
are computed with commonly used routines: i.e. 
the third order total variation diminishing Runge-Kutta scheme (TVD RK3) \cite{Shu-Osher-CL-1:1988}  for time discretization, and Godunov Hamiltonian~\cite{Rouy:1992} for the eikonal terms $\pm|\nabla\phi|$ with the fifth order WENO discretization \cite{jiang2000weighted} approximating $\nabla\phi$. We refer the readers to the book~\cite{LevelSet_OsherFedkiw} and~\cite{TO:Acta-Numerica:2005} for more detailed discussions and references.
We have also arranged our codes to be openly available on GitHUB.\footnote{ \texttt{https://github.com/GaZ3ll3/ibim-levelset}}

{A quick remark is in order regarding the algorithms used to generate an implicit representation of a 
"solvent excluded surface". 
Of course there are other approaches to generate the surfaces under the level set framework. While the general ideas appear to be similar, they are different in many details that could potentially influence performance of the algorithm that uses the prepared surfaces. We point out here that the procedure described in~\cite{MiThHeBoGi-CiCP13} is different from the one that proposed in Section~\ref{sec:vdW-F}. In particular, our method applies reinitialization \emph{after} removing the cavities, and does not need additional numerical solution of a Dirichlet problem for Laplace equation on irregular domain as in~\cite{MiThHeBoGi-CiCP13}. The reinitialization after cavity removal is essential to our proposed approach which requires distance function in the narrowband surrounding the surface. If there is no cavity, our procure does not require the reinitialization step.
We refer to~\cite{MiThHeBoGi-CiCP13,pan2009model}  for a more extensive review of other related algorithms. 
}

\subsection{Creating a signed distance function to the solvent accessible surface}
\label{sec:vdW-F}

From molecules description the van der Waals surface, $\Gamma^{vdW}$, is defined as the zero level set of 
\begin{equation}\label{vdW-levelset-fn}
	F(\bx) = \inf_{k} \left(\|\bx - \bz_k\| - r_k\right),
\end{equation}
where $\bz_k$ and $r_k$ denote respectively the coordinates of the molecule centers and their radii. 

From the van der Waals surface, we shall define the so-called solvent excluded surface, $\Gamma$, as the zero level set of a continuous function $\phi_{SES}$. 
$\phi_{SES}$ is computed by a simple inward eikonal flow, starting from an initial condition involving $F$, and is followed by a few iterations of the standard level set reinitialization steps. 
 See Figure~\ref{fig:surfs} for an illustration of this procedure in two dimensions.
 The details are described in the following subsections.
	
\subsubsection{Inward eikonal flow}
\label{sec:inward-eikonal-flows}

The van der Waals surface is extended outwards for a radius $\rho_0$ to {define} the so-called ``solvent accessible surface", which can be conveniently defined as the zero level set of $\phi_{SAS}$:
\begin{equation}
	\phi_{SAS}(\bx)= F(\bx) - \rho_0.
\end{equation}
The inward eikonal flow will produce a surface with smoothed out the corners when compared to the original van der Waals surface, while keeping most of its smooth parts unchanged.
For $0<t\le\rho_0$, we solve the following equation:
\begin{equation}
\begin{aligned}
	\frac{\partial {\tilde\phi}_{SAS}(\bx, t)}{\partial t} - \|\nabla \tilde\phi_{SAS}\| &= 0,~~~\bx\in U,\\
	\tilde\phi_{SAS}(\bx, 0) &= {\phi}_{SAS}(\bx),
\end{aligned}
\end{equation}
with zero Neumann boundary conditions.

\subsubsection{Cavities removal}
\label{sec:cavity-removal}
\begin{figure}[!htb]
\centering
\includegraphics[scale=0.25]{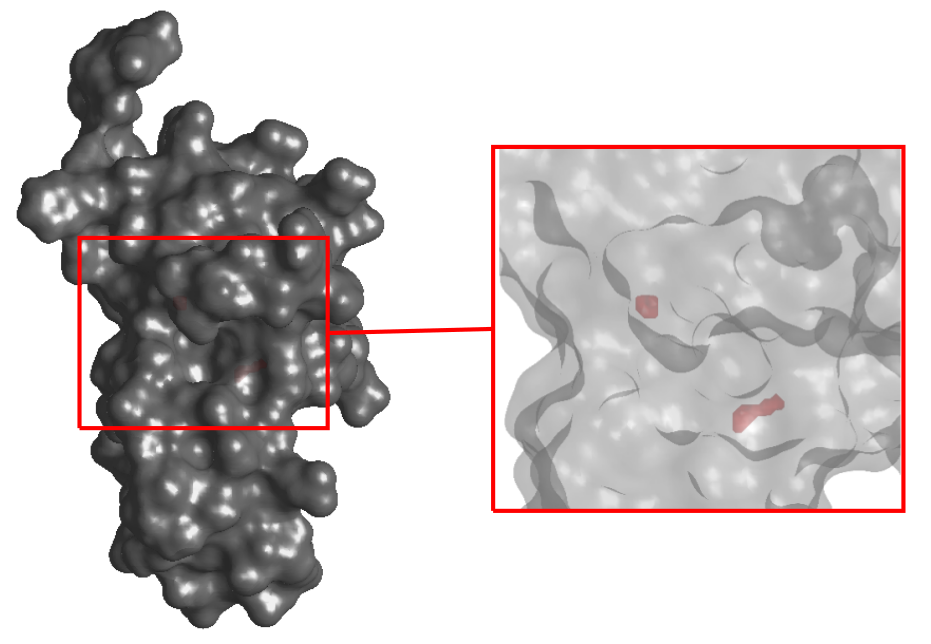}
\caption{Cavity in protein {\tt 1A63}. The gray surface is the ``solvent excluded surface'' used for computing the electric potential for protein {\tt 1A63}. In the right subfigure, the rendering of the gray surface is made semi-transparent in order to reveal the enclosed cavity surfaces (red). The regions enclosed by the red surfaces are the cavities to be removed.}
\label{fig:cavity}
\end{figure}

The zero level set of $\tilde\phi_{SAS}$ may contain some pieces of surfaces that isolate cavities that are believed to be void of solvent. Figure~\ref{fig:cavity} provides an example of such cavities in a protein that we used for computation. 
The cavity removal step uses a simple sweep to remove (the boundaries of) these regions and create a level set function, $\phi_{SES}$, that describes only the exterior, closed and connected surface --- the solvent excluded surface:
\[
\tilde\phi_{SAS}(\bx,\rho_0) \longrightarrow {\phi}_{SES}(\bx). 
\]

The cavity removal consists of following steps:
\begin{enumerate}
\item Identify a region $C_\epsilon$ that contains the cavities. $C_\epsilon$ is a superset of the cavities, containing points outside of the cavities that are within $\epsilon$ distance to the cavity surface. This can be done by a ``peeling" process: by moving the set of markers initially placed on the boundary of the computational domain inwards, using a breadth-first search (BFS) algorithm. The first layer of the surfaces defined by the zero level set of $\tilde\phi_{SAS}$ is defined to be the ``solvent excluded surface". We could therefore remove the remaining portion of $\tilde\phi_{SAS}$'s zero level sets, which are regarded {as} corresponding to the cavities.  From this process, one can easily compute a characteristic function supported on $C_\epsilon$. 
\item Remove the cavity region by modifying the values of $\tilde\phi_{SAS}$ in $C_\epsilon$:
\[
\phi_{SES}(\bx):=
\begin{cases}
\tilde{\phi}_{SAS}(\bx), & \bx\notin C_{\epsilon},\\
-\epsilon, & \bx\in C_{\epsilon}.
\end{cases}
\]
\end{enumerate}
	
\subsubsection{Reinitialization}
\label{sec:reinit}

The kinks on the solvent accessible surface (SAS) will lower the accuracy in the computation for 
$\phi_{SES}(\bx,\rho_0)$. In addition, the cavity removal step may introduce small jump discontinuities near the removed cavity. 
We  perform several iterations of the standard level set reinitialization to improve the 
equivalence of the computed $\phi_{SES}(\bx,\rho_0)$ and the signed distance function to $\Gamma$ (which $\phi_{SES}$ is supposed to be).
The reinitialization equation, first appeared in \cite{sussman1994level}, is defined as
\begin{equation}\label{eq:reinit}
\begin{aligned}
	\frac{\partial \tilde{\phi}_{SES}(\bx, t)}{\partial t} + \sgn_h(\phi_{SES}(\cdot,\rho_0))(|\nabla\tilde\phi_{SES}| - 1) &= 0,\\
	 \tilde{\phi}_{SES} (x,0) &= \phi_{SES}(x,\rho_0),
	\end{aligned}
\end{equation}
where the smoothed-out signum function is defined as
\begin{equation}\label{eq:signum}
	\sgn_h(\phi) = \frac{\phi}{\sqrt{\phi^2 + h^2}}.
\end{equation}
Suppose that the reinitialization equation is {solved} until $t=t_n$, i.e. $\tilde\phi_{SES}(\bx,t_n)$ is our approximation to the signed distance function $d_\Gamma(\bx)$, 
we shall {compute} $\Gamma_\epsilon$ by 
\begin{equation}
	\Gamma_\epsilon := \{ \bx\in\mathbb{R}^d: -\epsilon<\tilde\phi_{SES}(\bx, t_n) <\epsilon\}.
\end{equation}
We use the standard fifth order WENO for space and third order TVD-RK scheme for time to compute the reinitialization. In general, lower order schemes will result in larger perturbation of the zero level surface, which is not supposed to be moved.

The smoothness of the signum function $\mathrm{sgn}_h$ may influence the efficiency and effectiveness of the reinitialization procedure.
In our simulations, with the regularized signum function defined in \eqref{eq:signum}, 
it suffices to solve \eqref{eq:reinit} for $O(\epsilon)$ amount of time, in a neighborhood close to the zero level set. 
With $\Delta t=k_0 h$, and $\epsilon=k_1h$,  $k_0,k_1>0$, we run constant number of time steps for reinitialization, independent of $h$.
Since the fifth order WENO approximation of $\nabla\phi$ uses central differencing with seven grid nodes along each grid lines, the minimal number of time steps needed to create the signed distance in $\Gamma_\epsilon$ is $(k_1+4)/k_0$.
We refer the readers to \cite{Redistance_ChengTsai} for some more detailed discussion on reinitialization of level set functions and (closest point) extension of functions from $\Gamma$ to $\Gamma_\epsilon$., and an alternative higher order algorithm.

\subsection{Projections and weights}
\label{sec:weight-computations}

We locate all grid {points} $\bx_i\in U_h$ satisfying that $|\phi_{SES}(\bx_i)| < \eps$ and compute projections $\bx_i^{\ast}\in\Gamma$ by
\begin{equation}
	\bx_i^{\ast} = \bx_i -\tilde\phi_{SES}(\bx_i,t_n)\nabla \tilde\phi_{SES}(\bx_i,t_n).
\end{equation}
$\nabla \tilde\phi_{SES}(\bx_i,t_n)$ can be approximated either by standard central differencing or by the fifth order  \texttt{WENO} routines.
More precisely, on each grid node for each Cartesian coordinate direction, \texttt{WENO} returns two approximations of $\nabla\tilde\psi$, say $p_{-}$ and $p_+$, which are generalizations of the standard forward and backward finite differences of $\psi$. In our numerical simulations, we use 
\[
\nabla \tilde\phi_{SES} \approx \frac{p_-+p_+}{2}.
\]

For weight function $\delta_{\eps}$, we adopt {the} following cosine function with vanishing first moment, 
\begin{equation}\label{eq:cosine-kernel}
	\delta_{\eps}(\eta) = 
	\begin{cases}
		\frac{1}{2\epsilon}\left(1 + \cos \frac{\eta\pi}{\eps}\right), |\eta|<\epsilon,\\
		0, |\eta|\ge0.
	\end{cases}
\end{equation}
For general smooth nonlinear integrands and $\epsilon\sim{o}(1)$ for $h\rightarrow 0$, the above weight function provides at most 
second order in $h$ convergence. 
Since the chosen $\delta_\epsilon$ is an even function of the distance to the surface, it has one vanishing moment. 
Therefore, the zeroth order (in distance to the surface) approximation of the Jacobian will lead to a zero order in $\epsilon$ error.
See \cite{kublik2016extrapolative} for more in depth analysis on the properties of different choices of $\delta_\epsilon$.
In the simulations reported in this paper, we set  $J(\bx) \equiv 1$.


\subsection{Fast linear solvers}
\label{sec:FMM}

Equations \eqref{eq:IBIM-PBE}-\eqref{eq:weights} in Section~\ref{sec:IBIM},  together with the regularization of the kernels described in Section~\ref{sec:regularization}, one arrives at the final linear system:
\begin{equation}\label{eq:final-lin-sys}
\mathbf{\Lambda p} + \mathbf{KW p} = \mathbf{g},
\end{equation}
with $\mathbf{p}$ denoting the vector containing both $\bar\psi(\bx_j)$ and $\bar\psi_n(\bx_j)$, 
\[
\mathbf{\Lambda} :=\frac{1}{2}\left(\begin{array}{cc}
\left(1 + \frac{\epsilon_E}{\epsilon_I}\right)\mathbf{I} & \mathbf{0}\\
\mathbf{0} & \left(1 + \frac{\epsilon_I}{\epsilon_E}\right)\mathbf{I}
\end{array}\right),
\] and
$\mathbf{W}$ is a diagonal matrix defined by the weights $\omega_j:=J(\by_j)\delta_\epsilon(d_\Gamma(\by_j))$ as defined in Section~\ref{sec:IBIM}.

We solve this system by a standard \texttt{GMRES} algorithm. In the \texttt{GMRES} algorithm, we use the black-box fast multipole method (\texttt{BBFMM})~\cite{fong2009black} to accelerate the multiplication of the operator $\mathbf{K}$ to any vector. In particular, the solution of the diagonal part of~\eqref{eq:final-lin-sys} is used as a preconditioner. This means that the \texttt{GMRES} algorithm starts with the particular initial condition:
\[
	\mathbf{p}^{(0)}:= (\mathbf{\Lambda}+\mathbf{D})^{-1} \mathbf{g}, 
\]
where 
\[
D:=\left(\begin{array}{cc}
\mathbf{0} & \mathbf{0}\\
\mathbf{0} & \mathbf{0}
\end{array}\right)
\]
{comes} from the regularization of the kernels.

\subsection{Computing surface area}
\label{sec:surface-area}

In our IBIM approach, 
the evaluation of the surface area of $\Gamma$ is computed by $S^h_{\Gamma_\epsilon}[f]$ defined in~\eqref{EQ:IBIM sum} with $f\equiv 1$. 
See \cite{kublik2016integration}.
\subsection{Computing the polarization energy}
\label{sec:polarization-energy-computation}

The polarization energy $\mathcal{G}_{pol}$ of the system is given by
\begin{equation}\label{EQ:G pol}
	\mathcal{G}_{pol} = \frac{1}{2}\sum_{k=1}^{N_c}q_k \psi_{rxn}(\bz_k) 
\end{equation}
where $\psi_{rxn}(\bz_k)$ is computed by evaluating the following boundary integral at the center of atom $k$, $\bz_k$:
\begin{equation}\nonumber
	\psi_{rxn}(\bz) = \int_{\Gamma} \left(\left(\frac{\epsilon_E}{\epsilon_I}\frac{\partial G_{\kappa}(\bz, \by)}{\partial n(\by)} - \frac{\partial G_0(\bz, \by)}{\partial n(\by)}\right)\psi(\by) + (G_0(\bz, \by) - G_{\kappa}(\bz, \by))\frac{\partial \psi}{\partial n}(\by) \right)ds(\by). 
\end{equation}

In our IBIM approach, 
evaluation of this integral is computed by $S^h_{\Gamma_\epsilon}[f(\bz,\cdot)]$ defined in~\eqref{EQ:IBIM sum} with 
\begin{equation}
	f(\bz,\by) := \left(\frac{\epsilon_E}{\epsilon_I}\frac{\partial G_{\kappa}(\bz, P_\Gamma \by)}{\partial n(\by)} - \frac{\partial G_0(\bz, P_\Gamma\by)}{\partial n(\by)}\right)\bar\psi(\by) + (G_0(\bz, P_\Gamma\by) - G_{\kappa}(\bz, P_\Gamma\by))\frac{\partial \bar\psi}{\partial n}(\by).
\end{equation} 
	
\section{Numerical experiments}
\label{SEC:Num Exp}

We now perform some numerical experiments using the computational algorithm we developed. In all the numerical simulations, we set the dielectric parameters $\epsilon_I = 1.0$, $\epsilon_E = 80.0$ and Debye-H\"{u}ckel constant $\kappa = 0.1257\angstrom^{-1}$. We use the following parameters for the implicit boundary integral method:

\begin{tabular}{ll}
	$h$ &  denotes the grid spacing in the uniform Cartesian grids,\\
	$\epsilon\equiv 2h$ & denotes the width for the narrowband $\Gamma_{\epsilon}$,\\
	$\tau=h~\mathrm{or}~h/2$ & denotes the regularization parameter used in $K_{11},K_{12},K_{21},K_{22}$.
\end{tabular}

We set the tolerance in the \texttt{GMRES} algorithm to be $10^{-5}$, and use 4th order Chebyshev polynomials in the \texttt{BBFMM} preconditioner to achieve tolerance $10^{-4}$ there. In general, smaller $\tau$ results in more accurate approximations, if the resulting linear systems can be solved successfully. However, it cannot be too small with respect to the grid spacing, otherwise the resulting linear system becomes badly conditioned, as the current regularization approach becomes ineffective.

Most of the numerical experiments are performed on a desktop with quad-core CPU at 3.40GHz, 16GB RAM. The computations involving more than one million unknowns are performed on an older Linux computer with similar cache memory but sufficient RAM; for convenience in comparison, the timings presented in the tables below for simulations performed on this computer are scaled according to the clock speed and processor differences between the two computers. 
We think that such ad hoc scaling of timing is reasonable for the size of computations and the machines involved. We put an * sign next to the scaled CPU timings in the Tables.

In Section~\ref{sec:compare-to-MSMS}, we first compare the surface areas computed by our method to the ones computed by a published algorithm. In all later subsections, we present simulations of our algorithms with molecules of different sizes. In certain examples, we compare our computational results to the available published data. Particularly, we perform simulations on more realistic benchmark macromolecules taken from the RCSB Protein Data Bank (PDB)~\cite{berman2002protein}, and add missing heavy atoms through software \texttt{PDB2PQR}~\cite{dolinsky2004pdb2pqr}. The atom charges and radius parameters that we will be using in our simulations are all generated through force field \texttt{CHARMM}~\cite{BrBrOlStSwKa-JCC83}. The number of atoms reported in each subsection below {corresponds} to the number in the respective pqr file of each molecule.

Here are some general remarks on the numerical simulations using our algorithm. 	
In the tables presented in this section, the columns titled ``D.O.F.'' show the total number of unknowns in the discrete systems and provides a basis to observe the rate of convergence for the computed solutions, surface areas and polarization energies.
In the figures presented in the this section, the electrostatics computed on two different grid resolutions are painted on top of the respective
computed solvent excluded surfaces. In practice, the most important information that one would like to extract from such computations is the locations of the extrema; see e.g.~\cite{jones2003using,weiner1982electrostatic}. As the figures show, the extrema of the potentials computed at the coarser grid resolutions are already at the ``correct" locations on the surfaces.

As we shall see from the numerical accuracy study in Section~\ref{sec:single-ion}, the foremost bottleneck of the proposed algorithm 
is the low order regularization for the singular integrals. However, regularization is essential, and smaller amount of it (smaller values of $\tau$) leads to systems which require more \texttt{GMRES} iterations unless the the grid spacing $h$ is sufficiently small.  See, for instance, the simulations presented in Tables~\ref{tb:2aid},~\ref{tb:1f15} and~\ref{tb:1a2k}.
Despite the regularization issue, the boundary integrals can be computed very accurately if wider 
$\epsilon$ (with respect to the grid spacing $h$)  and the full expression of the Jacobian $J$ are used. However, wider $\epsilon$ implies a larger dense linear system {needs} to be solved. Most of the reported computation times are spent on the evaluation of the matrix-vector multiplications. Thus, in the simulations presented below, we choose a regime in which $\epsilon$ is narrow but sufficient in practice for the adopted simple quadrature to resolve the surface geometry. 
Finally, regarding to how small  $h$ should be for a given molecule and probe size $\rho_0$: 
	\[
	h\approx \min_{j=1,\cdots,N}\{\rho_0, r_j, 2\epsilon\}/7<\triangle
	\] 
where $\triangle$ is the minimal distance between ``different parts of the surface'' (think of the thin part of a dumb bell) and we define it as
\begin{equation}
\triangle = \inf_{\bx\in\Gamma} \{-d_{\Gamma}(\bx - s\bn(\bx)) | \ell(\bx) > s > 0, \text{ where }\bx -\ell(\bx)\bn(\bx)\in\Gamma \}.
\end{equation}
We shall see in the following examples, that our algorithm seems to perform well even {when} the discretized system is slightly outside of the above regime.	

%
%
%
%

\subsection{Molecular SES surface area}
\label{sec:compare-to-MSMS}

We compare the performance of our algorithm for calculating surface areas of different proteins with that of the \texttt{MSMS} (Michel Sanner's Molecular Surface) algorithm developed in~\cite{sanner1996reduced}. For \texttt{MSMS} algorithm, the probe radius is set to be $\rho_0=1.4\angstrom$ and the density parameter {is $1.0$} for mesh generation. We use the online implementation by the High-Performance Computing at the NIH group~\cite{StrucTools} to produce the data for \texttt{MSMS}. In Table~\ref{tb:ses-area}, we compare results from our method to these from \texttt{MSMS} for seven different proteins on a grid of size $128^3$. We observe that the surface areas computed by our algorithm are quite close to the \texttt{MSMS}'s approximate values in general. Since the \texttt{MSMS} results are only approximations to the true values, we did not attempt to tune algorithmic parameters (e.g. grid size, weight function) of our method to obtain results that are even closer to the \texttt{MSMS} results.

\begin{table}[htb!]
	\centering			
	\caption{Comparison of {surface} areas computed by the proposed  IBIM method and by \texttt{MSMS} for seven different proteins from the RCSB Protein Data Bank. }
	\medskip
	\begin{tabular}{|c c c|}
		\hline
		Protein id & Surface area  (IBIM)  & Surface area (\texttt{MSMS}) \\\hline
		\texttt{4INS}	& 4732    & 4761\\
		\texttt{1HJE}	& 825     & 801 \\
		\texttt{1A2B}    & 7540    & 7936\\
		\texttt{1PPE}    & 7979    & 8340 \\
		\texttt{2AID}    & 8061    & 8304\\
		\texttt{1F15}    & 22000   & 22725  \\
		\texttt{1A63}    & 6583    & 6659\\
		\hline
	\end{tabular}
	\label{tb:ses-area}
\end{table}


\subsection{The single ion model}\label{sec:single-ion}

We start with the single ion model developed in~\cite{kirkwood1934theory} to benchmark the solution accuracy of our numerical algorithm. We use three different relative errors, between exact and numerically represented quantities, to measure the quality of numerical solutions. They are defined as:
\begin{equation}
\begin{aligned}
\text{solution error} &= \sqrt{\frac{\int_{\Gamma} |\psi(\bx) - \psi^{\ast}(\bx)|^2  + |\frac{\partial\psi(\bx)}{\partial n} - \frac{\partial \psi^{\ast}(\bx)}{\partial n}|^2}{\int_{\Gamma} |\psi^{\ast}(\bx)|^2 + |\frac{\partial \psi^{\ast}(\bx)}{\partial n}|^2}},\\
\text{surface area error} &= \frac{|{\mathcal{A}} - {\mathcal{A}}^{\ast}|}{{\mathcal{A}}^{\ast}},\\
\text{energy error} &= \frac{|\mathcal{G}_{pol} - {\mathcal{G}}^{\ast}_{pol}|}{{\mathcal{G}}^{\ast}_{pol}}.
\end{aligned}
\end{equation}

For a single atom with radius $r$ and charge $q$, the solution to the Poisson-Boltzmann equation is given as~\cite{kirkwood1934theory}
\begin{equation}
\psi^{\ast}(\bx) = 
\begin{cases}
\dfrac{q}{4\pi \epsilon_I |\bx|} + \dfrac{q}{4\pi r}\left(\dfrac{1}{\epsilon_E(1+\kappa r)} - \dfrac{1}{\epsilon_I}\right),& \text{if } |\bx| < r\\
\dfrac{qe^{-\kappa (|\bx| - r)}}{4\pi \epsilon (1+\kappa r)|\bx|},              & \text{otherwise}
\end{cases}
\end{equation}
We can therefore compute the associated polarization energy
\begin{equation}
{\mathcal{G}}^{\ast}_{pol} = \dfrac{q^2}{8\pi r} \left(\frac{1}{\epsilon_E(1+\kappa r)} - \frac{1}{\epsilon_I}\right),
\end{equation}
using the fact that the surface is a sphere with area ${\mathcal{A}}^{\ast} = 4\pi r^2$. We set the atom's radius to be $r=1\angstrom$ and assigned charge to be $q=1e_c$.

\begin{table}[htb!]
	\caption{Benchmarking errors in solution of the single ion model.  }
	\medskip
	\centering
	\begin{tabular}{|cccrcccc|}
		\hline\rule{0pt}{2.5ex}
		grid size & $h (\angstrom)$  & $\tau/h$ & D.O.F.  & \texttt{GMRES} & solution error & area error & energy error \\\hline
		\rule{0pt}{2.2ex}$64^3$ & 3.91E$-$1 & $1$& 22,756 & 4  & 1.11E$-$02 &1.24E$-$03 & 1.21E$-$02 \\
		\rule{0pt}{1ex}$128^3$ &1.95E$-$1 & $1$  & 91,564  & 3  & 6.90E$-$03 &2.88E$-$04 & 5.68E$-$03 \\
		\rule{0pt}{1ex}$256^3$  &9.77E$-$2 & $1$ & 366,868 & 3 & 3.57E$-$03 & 5.01E$-$05 & 2.90E$-$03\\
		\hline
		\rule{0pt}{2.2ex}$64^3$  &3.91E$-$1  & $0.5$ & 22,756 & 2 & 5.93E$-$03 & 1.24E$-$03 & 5.98E$-$03 \\
		\rule{0pt}{1ex}$128^3$  & 1.95E$-$1& $0.5$  & 91,564  & 4  & 3.45E$-$03 &2.88E$-$04 & 2.49E$-$03 \\
		\rule{0pt}{1ex}$256^3$  & 9.77E$-$2& $0.5$ & 366,868 & 3  & 2.15E$-$03 & 5.01E$-$05 & 1.31E$-$03\\
		\hline
	\end{tabular}
	\label{tb:single-ion-model}
\end{table}
We performed simulations under different mesh and IBIM parameters. The results are summarized in Table~\ref{tb:single-ion-model}. Our method converges in very small numbers (usually $3\sim 4$) of iterations. This benchmark calculation shows that our numerical algorithms can indeed achieve similar solution accuracies to those achieved by other algorithms developed recently~\cite{BaChRa-SIAM11,BoFeZh-JPCB02}.

\subsection{Protein \texttt{1A63}}

In this numerical example, we compute the polarization energy for protein \texttt{1A63}, the E.Coli Rho factor, of the Protein Data Bank. The protein has $2065$ atoms with different radii. The information on the locations and radii of the atoms are all available in~\cite{berman2002protein}. 

\begin{figure}[!htb]
	\centering
	\includegraphics[width=0.45\linewidth, height=0.36\linewidth]{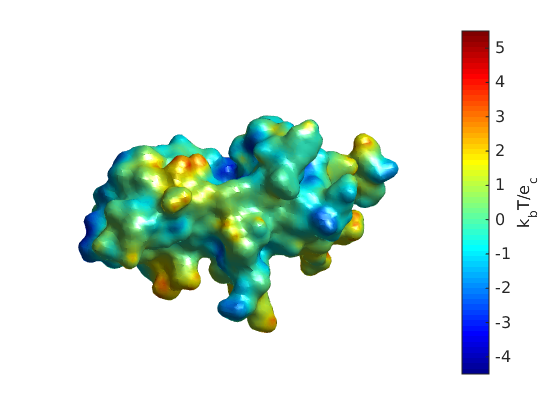}
	\includegraphics[width=0.45\linewidth, height=0.36\linewidth]{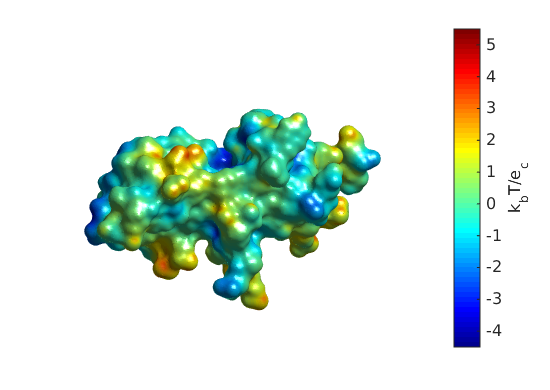}
	\caption{The electrostatic potential on the surface of the PDB-\texttt{1A63} protein. Left: on grid $128^3$. Right: on grid $256^3$. }
	\label{fig:1A63}
\end{figure}
In Figure~\ref{fig:1A63} we plot the  potential $\bar\psi$ on the constructed ``solvent excluded surface'', computed on two different grids, $128^3$ (left) and $256^3$ (right). Further computational results are tabulated in Table~\ref{tb:1a63}. The computed values of the polarization energy $\cG_{pol}$ can be compared to the existing estimations, $\mathcal{G}_{pol}^{\mathrm{TABI}} = -2374.64~\text{kcal/mol}$ from the treecode-based boundary integral solver \texttt{TABI}~\cite{GeKr-JCP13} and  $\mathcal{G}_{pol}^{\mathrm{APBS}} = -2350.58~\text{kcal/mol}$ from the finite difference solver \texttt{APBS}~\cite{unni2011web}. The computational results are tabulated in Table~\ref{tb:1a63}. 

\begin{table}[htb!]
	\centering
	\caption{Numerical results on protein \texttt{1A63} under different algorithmic parameters. }
	\medskip
	\begin{tabular}{|cccrcccc|}
		\hline\rule{0pt}{2.5ex}
		grid size & $h(\angstrom)$& $\tau/h$& D.O.F.  & \texttt{GMRES}  & $\mathcal{G}_{pol}$ (kcal/mol)  & CPU (s) & area $(\angstrom^2)$ \\\hline
		\rule{0pt}{2.2ex}$128^3$& 6.11E$-$1 &1 & 71,597  & 12          & -2392.22 & 606.8 & 6583\\
		\rule{0pt}{1ex}$256^3$ & 3.05E$-$1& 1 & 293,627  & 13          & -2366.40 & 3041 & 6801\\
		\hline
		\rule{0pt}{2.2ex}$128^3$& 6.11E$-$1&0.5 & 71,597  & 13         & -2345.41 &  772.3 & 6583\\
		\rule{0pt}{1ex}$256^3$ & 3.05E$-$1& 0.5 & 293,627  & 14          & -2347.74 &  3808 &6801\\
		\hline
	\end{tabular}
	\label{tb:1a63}
\end{table}

\subsection{Protein \texttt{2AID}}

\begin{figure}[!htb]
	\centering
	\includegraphics[width=0.45\linewidth, height=0.36\linewidth]{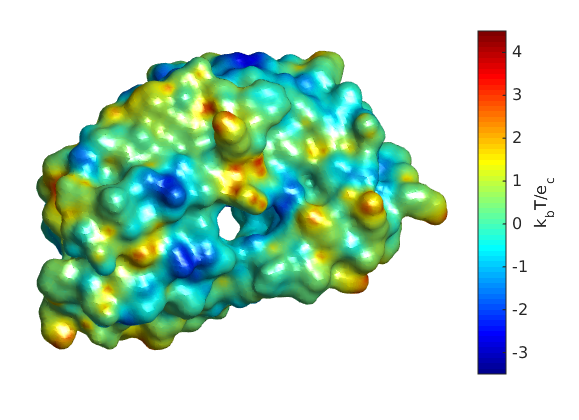}
	\includegraphics[width=0.45\linewidth, height=0.36\linewidth]{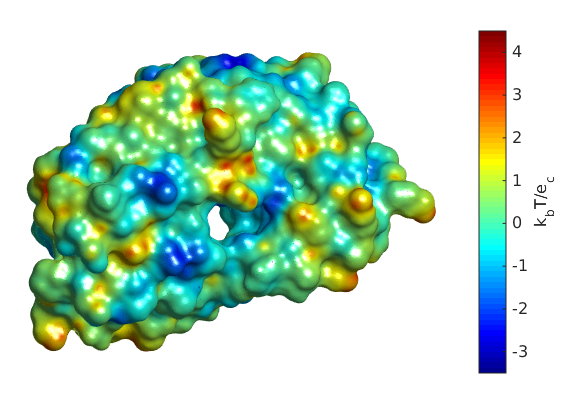}
	\caption{The electrostatic potential on the surface of protein \texttt{2AID} for two grids of size $128^3$ (left) and $256^3$ (right).}
	\label{fig:2aid}
\end{figure}

Here we compute the polarization energy for protein \texttt{2AID}, a non-peptide inhibitor complexed with HIV-1 protease. This protein has $3130$ atoms. In Figure~\ref{fig:2aid} we plot the  potential $\bar\psi$ on the constructed ``solvent excluded surface'' of this protein, computed on two different grids. Further computational results are tabulated in Table~\ref{tb:2aid}. 

\begin{table}[htb!]
	\caption{Numerical results on protein \texttt{2AID} under different algorithmic parameters.}
	\medskip
	\centering
	\begin{tabular}{|cccrcccc|}
		\hline\rule{0pt}{2.5ex}
		grid size &$h(\angstrom)$& $\tau/h$& D.O.F.  & \texttt{GMRES} & $\mathcal{G}_{pol}$ (kcal/mol)  & CPU (s) & area $(\angstrom^2)$ \\\hline
		\rule{0pt}{2.2ex}$128^3$& 5.80E$-$1&1 & 97,108  & 13  & -2318.69  & 940.7 & 8061\\
		\rule{0pt}{1ex}$256^3$ &2.90E$-$1& 1 & 397,930 & 14  & -2321.72 & 4521& 8335\\
		\hline
		\rule{0pt}{2.2ex}$128^3$& 5.80E$-$1 &0.5 & 97,108  & 24  & -2282.14  & 1745  &8601\\
		\rule{0pt}{1ex}$256^3$ &2.90E$-$1& 0.5 & 397,930 & 15  & -2306.70 & 4906&8335\\
		\hline
	\end{tabular}
	\label{tb:2aid}
\end{table}

\subsection{Protein \texttt{1F15}} 

\begin{figure}[!htb]
	\centering
	\includegraphics[width=0.45\linewidth, height=0.36\linewidth]{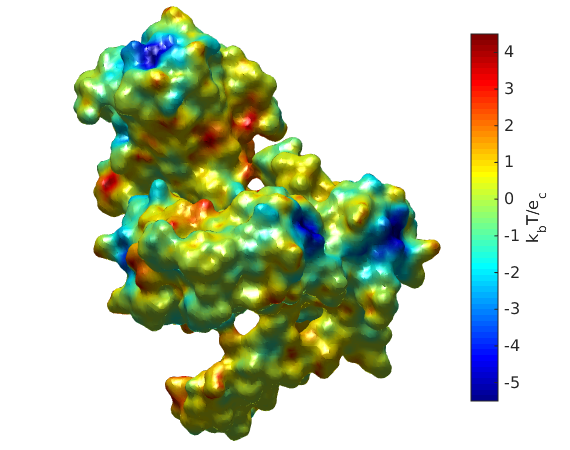}
	\includegraphics[width=0.45\linewidth, height=0.36\linewidth]{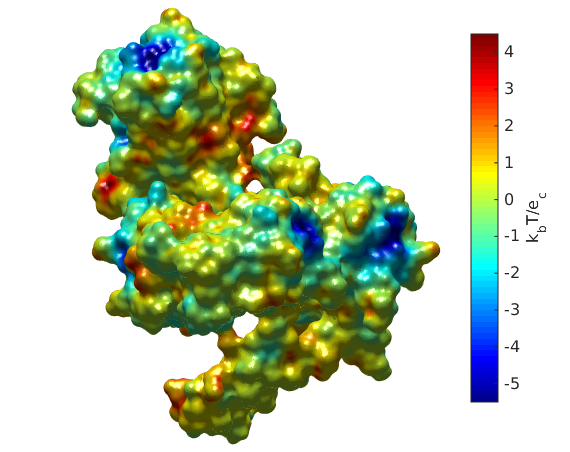}
	\caption{The electrostatic potential on the surface of protein \texttt{1F15} for two different grids of sizes $128^3$ (left) and $256^3$ (right).}
	\label{fig:1f15}
\end{figure}

In this example, we compute the polarization energy for protein~\texttt{1F15}, the cucumber mosaic virus. The protein has $8494$ atoms. 
In Figure~\ref{fig:1f15} we plot the potential $\bar\psi$ on the constructed ``solvent excluded surface'', computed on two different grids. Further computational results are tabulated in Table~\ref{tb:1f15}.

\begin{table}[htb!]
	\caption{Numerical results on protein \texttt{1F15} under different algorithmic parameters.}
	\medskip
	\centering
	\begin{tabular}{|cccrcccc|}
		\hline\rule{0pt}{2.5ex}
		grid size &$h(\angstrom)$& $\tau/h$& D.O.F.  & \texttt{GMRES} & $\mathcal{G}_{pol}$ (kcal/mol)  & CPU (s) & area $(\angstrom^2)$ \\\hline
		\rule{0pt}{2.2ex}$128^3$& 7.72E$-$1&1& 147,463  & 17 & -7770.00  & 2586 & 22000 \\
		\rule{0pt}{1ex}$256^3$ &3.86E$-$1& 1& 613,726  & 21 & -7818.83 & 12357 & 22847\\
		\rule{0pt}{1ex}$512^3$ &1.93E$-$1 & 1 & 2,497,309 & 25 & -7891.05 & 45681*& 23238\\
		\hline
		\rule{0pt}{2.2ex}$128^3$ &7.72E$-$1&0.5& 147,463  & 31 & -7682.67  & 4667 &22000\\
		\rule{0pt}{1ex}$256^3$ &3.86E$-$1& 0.5& 613,726  & 26 & -7774.76 & 14129&22847\\
		\rule{0pt}{1ex}$512^3$ &1.93E$-$1 & 0.5 & 2,497,309 & 29 & -7875.86 & 52687*& 23238\\
		\hline
	\end{tabular}
	\label{tb:1f15}
\end{table}
\subsection{Protein \texttt{1A2K}}

In this numerical example, we compute the polarization energy for protein \texttt{1A2K}, the GTPase RAN-NTF2 complex. The protein has $13627$ atoms. 

\begin{figure}[!htb]
	\centering
	\includegraphics[width=0.45\linewidth, height=0.36\linewidth]{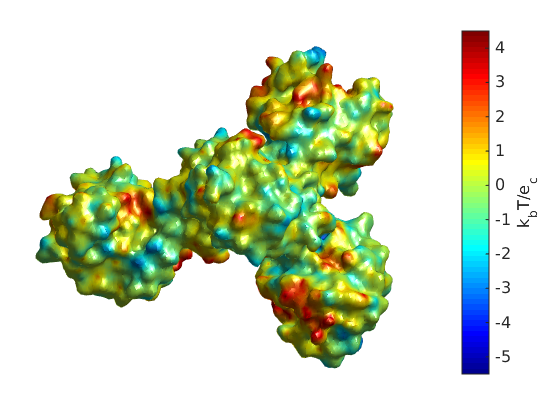}
	\includegraphics[width=0.45\linewidth, height=0.36\linewidth]{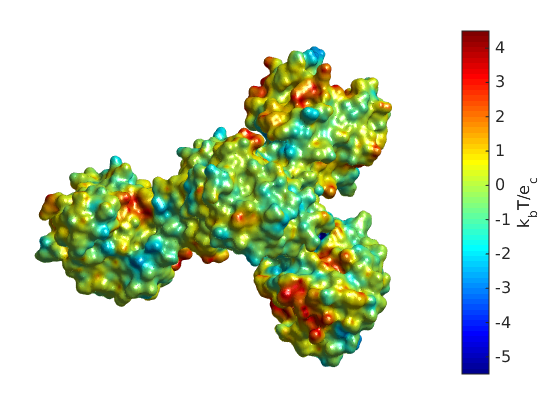}
	\caption{The electrostatic potential on the surface of protein \texttt{1A2K} for two different grids of sizes  $128^3$ (left) and $256^3$ (right).}
	\label{fig:1a2k}
\end{figure}

In Figure~\ref{fig:1a2k} we plot the  potential $\bar\psi$ on the constructed ``solvent excluded surface'', computed on two different grids. Further computational results are tabulated in Table~\ref{tb:1a2k}.

\begin{table}[htb!]
	\caption{Numerical results on protein \texttt{1A2K} under different algorithmic parameters.}
	\medskip
	\centering
	\begin{tabular}{|cccrcccc|}
		\hline\rule{0pt}{2.5ex}
		grid size & $h(\angstrom)$ & $\tau/h$& D.O.F.  & \texttt{GMRES} & $\mathcal{G}_{pol}$ (kcal/mol)  & CPU (s) & area $(\angstrom^2)$ \\\hline
		\rule{0pt}{2.2ex}$128^3$&1.08E$+$0 &1& 99,783  & 12 & -7190.38  & 1032 & 29497\\
		\rule{0pt}{1ex}$256^3$ &5.42E$-$1& 1 & 429,451  & 17 & -8902.62 & 6120 & 31501\\
		\rule{0pt}{1ex}$384^3$ & 3.61E$-$1& 1 & 983,418 & 17 & -8920.93  & 11948 & 32047 \\
		\rule{0pt}{1ex}$512^3$ &2.71E$-$1& 1& 1,765,673 & 18  &  -8963.12 & 25496* & 32364\\
		\hline
		\rule{0pt}{2.2ex}$128^3$&1.08E$+$0 &0.5 & 99,783  & 21  & -9004.68  & 2309 & 29497\\
		\rule{0pt}{1ex}$256^3$ &5.42E$-$1& 0.5 & 429,451  & 43  & -8789.45 &  15528 &31501\\
		\rule{0pt}{1ex}$384^3$ &3.61E$-$1& 0.5& 983,418   & 40  & -8859.47 &  26030 &32047\\
		\rule{0pt}{1ex}$512^3$ &2.71E$-$1& 0.5& 1,765,673 & 23  &   -8921.64 &  33558* & 32364\\
		\hline
	\end{tabular}
	\label{tb:1a2k}
\end{table}

\subsection{Protein: \texttt{1PMA}}
\begin{figure}[!htb]
	\begin{center}
		\includegraphics[width=0.45\linewidth, height=0.36\linewidth]{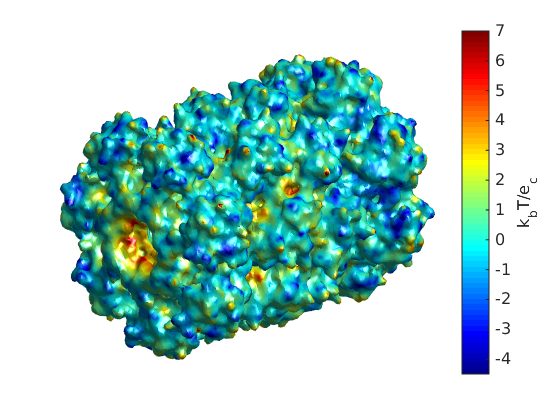}
		\includegraphics[width=0.45\linewidth, height=0.36\linewidth]{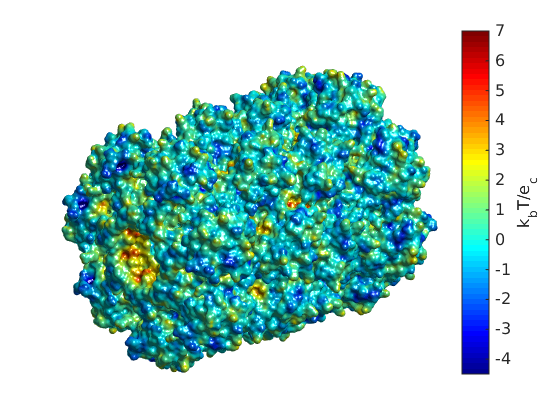}
	\end{center}
	\caption{The electrostatic potential on molecular surface for the proteasome from thermoplasma acidophilum (PDB id: 1PMA). Left: the potential computed on a $128^3$ grid. Right: the potential computed on a $512^3$ grid. }\label{fig:1pma}
\end{figure}

In this example, we compute the polarization energy for proteasome from thermoplasma acidophilum (PDB id: 1PMA) with 93017 atoms. 
In Figure~\ref{fig:1pma} we plot the  potential $\bar\psi$ on the constructed ``solvent excluded surface'', computed on two different grids. Further computational results are tabulated in Table~\ref{tb:1pma}. 

\begin{table}[htb!]
	\caption{Protein 1PMA. IBIM's result of relative error in polarization energy and total run time w.r.t  different grid sizes.}\bigskip
	\centering
	\begin{tabular}{|cccrcccr|}
		\hline\rule{0pt}{2.5ex}
		grid size &$h(\angstrom)$& $\tau/h$& D.O.F.  & \texttt{GMRES}   & $\mathcal{G}_{pol}$ (kcal/mol)  & CPU (s) & area $(\angstrom^2)$ \\\hline
		\rule{0pt}{2.2ex}$128^3$& 1.67E$+$0&1 & 222,478  & 15  & -15544.35  & 3360 & 1.6014E$+$5\\
		\rule{0pt}{1ex}$256^3$ &8.35E$-$1& 1 & 1,026,938 & 18  &  -46865.78 & 14141* & 1.8034E$+$5\\
		\rule{0pt}{1ex}$384^3$ &5.56E$-$1& 1 & 2,429,367 & 20  &  -50071.90 & 42927* & 1.8851E$+$5\\
		\rule{0pt}{1ex}$512^3$ & 4.17E$-$1& 1 & 4,418,314 & 22 &  -50144.63 & 70880* & 1.9262E$+$5\\
		\hline
	\end{tabular}
	\label{tb:1pma}
\end{table}

\section{Concluding remarks}
\label{SEC:Concl}

We present in this paper a new numerical method for solving the boundary value problem of the linearized Poisson-Boltzmann equation, which is widely used to model the electric potential for macromolecules in solvent. Our new method relies on the standard level set method~\cite{LevelSet_OsherFedkiw,osher_sethian88} for preparing the distance function to the ``molecular surface''. Contrary to the typical level set method, in which some partial differential equations are discretized with some suitable boundary conditions, ours {involves} the solution of an integral equation which is derived from an implicit boundary integral formulation~\cite{KTT}. Similar to the typical level set methods, and contrary to the typical boundary integral methods, the proposed method {involves} computation only with functions defined on uniform Cartesian grids. Our numerical simulations show that in addition to the flexibility that comes from the level set methods, the proposed method can be as computationally efficient as other boundary integral based algorithms. We show by our numerical simulations that the solutions of the resulting linear systems can be accelerated easily by some existing fast multipole methods. Furthermore, the eikonal flow and reinitialization in Stage (1) of the proposed algorithm rely on widely available explicit solvers and can be trivially parallelized. 

There are several possible improvements that could be investigated in the future. First of all, the quadrature for the implicit boundary integral formulation can be improved to increase the order of accuracy. This includes improvement of the regularization of the kernel singularities and the use of full expression of the Jacobian $J$. One may also consider different grid geometries, as the underlying mathematical formulation {does} not require uniform Cartesian grids. For example, the adaptive oct-tree structure used in~\cite{HeGi-JCP11} or radial basis functions may be considered.

As all the presented simulations were computed on two moderate desktop computers, 
the reported results show the potential of the proposed method for molecular dynamics simulations involving very large molecules.

Finally, let us mention that the numerical method we proposed here can be generalized to solve many similar model problems for electrostatics in related areas of electrochemistry~\cite{HoLiLiEi-JPCB12,MaPeAg-JCP88,MaMu-IJHMT09,HaGaRe-JCP16,HeGaLeRe-SIAM15}.

\section*{Acknowledgments}

Zhong and Ren are partially supported by the National Science Foundation through grants DMS-1321018 and DMS-1620473. Tsai is partially supported by the National Science Foundation through grants DMS-1318975 and DMS-1720171, and Army Research Office Grant No. W911NF-12-1-0519.
We would like to thank Professor Chandrajit Bajaj for motivating conversations on topics related to this paper.

{
\bibliography{IBIM-LS}
\bibliographystyle{siam}
}

\end{document}